\documentclass[11pt,a4paper]{article}
\usepackage{amsmath}
\usepackage{amsfonts}
\usepackage{amssymb}
\usepackage{amsthm}
\usepackage{mathrsfs}
\usepackage{hyperref}
\usepackage{url}
\usepackage[top=2cm, bottom=2cm, left=2.5cm, right=2.5cm]{geometry}
\newtheorem{thm}{Theorem}[section]
\newtheorem{lem}[thm]{Lemma}

\newtheorem{cor}[thm]{Corollary}

\theoremstyle{definition}

\theoremstyle{remark}
\newtheorem{rem}{Remark}

\newcommand\Tr{\operatorname{Tr}}

\newcommand\Leb{\operatorname{Leb}}

\usepackage{tikz}
\usetikzlibrary{trees}
\usepackage{url}

\begin{document}
\pgfsetxvec{\pgfpoint{8pt}{0}}
\pgfsetyvec{\pgfpoint{0}{8pt}}

\title{A convergence result on the lengths of Markovian loops}
\author{Yinshan Chang\footnote{Max Planck Institute for Mathematics in the Sciences, Inselstr. 04103, Leipzig, ychang@mis.mpg.de}}
\date{}
\maketitle

\begin{abstract}
Consider a sequence of Poisson point processes of non-trivial loops with certain intensity measures $(\mu^{(n)})_n$, where each $\mu^{(n)}$ is explicitly determined by transition probabilities $p^{(n)}$ of a random walk on a finite state space $V^{(n)}$ together with an additional killing parameter $c^{(n)}=e^{-a\cdot\sharp V^{(n)}+o(\sharp V^{(n)})}$. We are interested in asymptotic behavior of typical loops. Under general assumptions, we study the asymptotics of the length of a loop sampled from the normalized intensity measure $\bar{\mu}^{(n)}$ as $n\rightarrow\infty$. A typical loop is small for $a=0$ and extremely large for $a=\infty$. For $a=(0,\infty)$, we observe both small and extremely large loops. We obtain explicit formulas for the asymptotics of the mass of intensity measures, the asymptotics of the proportion of big loops, limit results on the number of vertices (with multiplicity) visited by a loop sampled from $\bar{\mu}^{(n)}$. We verify our general assumptions for random walk loop soups on discrete tori and truncated regular trees. Finally, we consider random walk loop soups on complete graphs. Here, our general assumptions are violated. In this case, we observe different asymptotic behavior of the length of a typical loop.
\end{abstract}

\section{Introduction}
 Poisson ensembles of Markovian loops were introduced informally by Symanzik \cite{Symanzik} and then by Lawler and Werner \cite{LawlerWernerMR2045953} for 2D Brownian motion. An extensive investigation of the loop soup on finite and infinite graphs was done by Le Jan \cite{LeJanMR2815763} for reversible Markov processes, and then by Sznitman \cite{SznitmanMR2932978} in relation with random interlacement (see \cite{SznitmanMR2680403} or \cite{DrewitzRathSapozhnikovMR3308116}).

 It is natural to study the typical behavior of loops in loop soups. The following questions were raised by Le Jan: for a loop sampled according to the normalized intensity measure $\bar{\mu}$, what can one say about its typical size? In particular, does the loop cover a positive proportion of the space? These questions are related to the length of the loop (i.e. the number of vertices with multiplicity visited by the loop). Indeed, by \cite{AleliunasKarpLiptonLovaszRackoffMR598110}, the cover time for the SRW on a general graph with the vertex set $V$ is $O(\sharp V^3)$, which implies that a randomly chosen Markovian loop visits all the vertices with high probability if its length exceeds certain power of $\sharp V$. This motivates our study on the lengths of loops.

 Given a sequence of Poisson point processes of loops on a sequence of increasing graphs, we are interested in the typical length of the loops. The intensity measures considered in the paper are given by Markov chains with uniform killing rates. When the killing rates go to zero super-exponentially, typical loops are extremely large. On the other hand, when the killing rates decrease to zero sub-exponentially fast, the length of a randomly chosen loop is tight. We are particularly interested in the intermediate case where two types of typical loops appear simultaneously: loops of bounded size and loops of size much bigger than the size of the state space of Markov chains.

 Following \cite{LeJanMR2815763}, we consider measures on the space of loops associated with these Markov chains. We first introduce the notation of based loops and loops. By a based loop $\dot{\ell}$, we mean an element $(x_1,\ldots,x_k)\in (V^{(n)})^k$. Two based loops are equivalent if they coincide up to a circular permutation, e.g. $(x_1,x_2,x_3,x_4)=(x_2,x_3,x_4,x_1)$. A loop is an equivalence class of based loops. The based loop functionals that we will be interested in will be loop functionals, i.e. they are invariant under circular permutations. For example, the \emph{length} $|\dot{\ell}|$ of a based loop $\dot{\ell}=(x_1,\ldots,x_k)$, which is defined by $|(x_1,\ldots,x_k)|=k$. We use the same notation for the length of a loop. We define the mass of $\dot{\ell}$ under the \emph{based loop measure} $\dot{\mu}^{(n)}$ by
\begin{equation}\label{eq: based loop measure}
\dot{\mu}^{(n)}(\dot{\ell})=\frac{1}{k}(1+c^{(n)})^{-k}\cdot p^{(n)}(x_1,x_2)p^{(n)}(x_2,x_3)\cdots p^{(n)}(x_k,x_1),
\end{equation}
where $(p^{(n)}(x,y))_{x,y}$ are the transition probabilities of an \emph{irreducible} Markov chain $X^{(n)}=(X^{(n)}_{k})_k$ on a finite state space $V^{(n)}$ and
\[c^{(n)}=e^{-a\cdot\sharp V^{(n)}+o(\sharp V^{(n)})}\quad (a>0)\]
plays the role of killing for the Markov chain $X^{(n)}$.

We normalize $\dot{\mu}^{(n)}$ and get a probability measure $\bar{\dot{\mu}}^{(n)}$. The corresponding push-forward measures on the space of loops are denoted by $\mu^{(n)}$ and $\bar{\mu}^{(n)}$. We refer to \cite{LeJanMR2815763} and \cite{SznitmanMR2932978} for more details on Markovian loop measures.

The main object of interest in this paper is the limit of $|\ell|$ under $\bar{\mu}^{(n)}$. Let \[\nu^{(n)}=\sum_{i=1}^{\sharp V^{(n)}}\delta_{\lambda_{n,i}}/\sharp V^{(n)}\]
be the \emph{empirical distribution of the eigenvalues of the matrix $(p^{(n)}(x,y))_{V^{(n)}}$}, where $\lambda_{n,i}$ are the eigenvalues of $p^{(n)}$ and $\delta_x$ is the Dirac mass at $x$ for $x\in \mathbb{C}$. By the Perron-Frobenius theorem, $\nu^{(n)}$ is supported in the unit disk. By tightness, there always exists a convergent subsequence of $\nu^{(n)}$. Therefore, we assume as \emph{hypothesis (H1)} that
\begin{equation}\label{eq: cvg assumption}
   \lim_{n\rightarrow\infty}\nu^{(n)}=\nu.
 \end{equation}
By the Perron-Frobenius theorem, there exists a unique stationary distribution $\pi^{(n)}$ for each $X^{(n)}$. Our second \emph{hypothesis (H2)} is that these $X^{(n)}$ are more or less similar to Markov chains on graphs of uniformly bounded degrees: there exist $n$-independent constants $c_1,c_2>0$ such that
 \begin{equation}
  \min_{x\in V^{(n)}}\pi^{(n)}(x)\cdot\sharp V^{(n)}\geq c_1,\quad\min_{\substack{x,y\in V^{(n)}\\p^{(n)}(x,y)+p^{(n)}(y,x)>0}}p^{(n)}(x,y)\geq c_2.
 \end{equation}
Our main result is a limit result on the typical size of a loop sampled from $\bar{\mu}^{(n)}$.
\begin{thm}\label{thm: thm1}
We suppose (H1) and (H2). Let $(q_n)_n$ be such that $\lim_{n\rightarrow\infty}q_n=\infty$ and that $\log q_n=o(\sharp V^{(n)})$ as $n$ increases to infinity. Then,
\begin{itemize}
\item[a)] $\int |\log(1-x)|\nu(\mathrm{d}x)\leq 3+8\sqrt{2}c_1^{-1}c_2^{-4}$,
\item[b)] $\lim_{n\rightarrow\infty}||\mu^{(n)}||/\sharp V^{(n)}=a+\int-\log(1-x)\nu(\mathrm{d}x)$ where $||\mu^{(n)}||$ is the total mass of $\mu^{(n)}$,
\item[c)] $\lim_{n\rightarrow\infty}\bar{\mu}^{(n)}(|\ell|\geq q_n)=a/(a+\int-\log(1-x)\nu(\mathrm{d}x))$,
\item[d)] $\lim_{n\rightarrow\infty}\bar{\mu}^{(n)}(|\ell|=j\big||\ell|\leq q_n)=j^{-1}\int x^j\nu(\mathrm{d}x)/\int-\log(1-x)\nu(\mathrm{d}x)$,
\item[e)] $\lim_{n\rightarrow\infty}\bar{\mu}^{(n)}(\log|\ell|/\sharp V^{(n)}\in \mathrm{d}x\big||\ell|\geq q_n)=1_{0<x<a}\cdot\frac{1}{a}\cdot\mathrm{d}x$.
\end{itemize}
\end{thm}
In Sections \ref{sec: torus} and \ref{sec: balls in regular trees}, we compute explicitly $\nu$ and $\int-\log(1-x)\nu(\mathrm{d}x)$ on tori and trees. The proof of Theorem~\ref{thm: thm1} is based on an upper bound on the transition functions of Markov processes and a direct calculation.

In contrast, when $X^{(n)}$ is the simple random walk on the complete graph of $n$ vertices. Theorem \ref{thm: thm1} is not applicable as (H2) fails. By an explicit calculation, we get
\begin{thm}\label{thm: thm2}
Suppose that $X^{(n)}$ is the simple random walk on a complete graph of $n$ vertices and $c^{(n)}=n^{-a+o(1)}$ for some $a>0$. Then, we have that $\lim_{n\rightarrow\infty}||\mu^{(n)}||/\log n=a$ and that $\lim_{n\rightarrow\infty}\bar{\mu}^{(n)}(\log|\ell|/\log n\in \mathrm{d}x)=1_{0<x<a}\cdot\frac{1}{a}\cdot\mathrm{d}x$.
\end{thm}
By a comparison with the coupon collector problem, we get
\begin{cor}\label{cor: cor of thm2}
 In the same setting as Theorem~\ref{thm: thm2}, let $\mathcal{L}_{n}$ be a Poisson point process of loops with intensity $\mu^{(n)}/\log n$. Then, $\sum_{\ell\in\mathcal{L}_n}1_{\ell\in\mathcal{C}}$ converges in distribution to a Poisson random variable with mean $\max(a-1,0)$, where $\mathcal{C}$ stands for the set of loops which cover all the vertices inside the graph.
\end{cor}

\emph{Organization of the paper}: Section~\ref{sec: proof of thm1} is devoted to prove Theorem~\ref{thm: thm1}. Then, we analyze two examples and calculate the limit measures $\nu$: random walk loop soups on discrete tori in Section~\ref{sec: torus} and random walk loop soups on balls in a regular tree in Section~\ref{sec: balls in regular trees}. We consider random walk loop soups on complete graphs and prove Theorem \ref{thm: thm2} and Corollary \ref{cor: cor of thm2} in the last section.

\section{Proof of Theorem~\ref{thm: thm1}}\label{sec: proof of thm1}
By (H2), the random walk $X^{(n)}$ is either aperiodic or $2$-periodic. For simplicity, we assume that these \emph{$X^{(n)}$ are all aperiodic}. The argument for periodic case is quite similar and is left to the reader. The key to prove Theorem~\ref{thm: thm1} is the heat kernel bounds in \cite{MorrisPeresMR2198701}, which is stated for lazy random walks\footnote{By lazy random walks, we mean random walks with the transition probabilities $((1-\gamma)p(x,y)+\gamma\cdot1_{x=y}))_{x,y}$, where $\gamma\in(0,1]$ and $(p(x,y))_{x,y}$ is the transition probabilities of some random walks.}. Under (H2), the laziness assumption is satisfied by $(X^{(n)}_{2k})_{k}$ and \cite[Theorem~4]{MorrisPeresMR2198701} gives an upper bound on $\int |x|^k\nu^{(n)}(\mathrm{d}x)$. This bound is then used to bound $\int |\log(1-x)|\nu(\mathrm{d}x)$ and to control the convergence rate of $\Tr (p^{(n)})^k$ towards $1$, see Lemma~\ref{lem: trace}. Then, Theorem~\ref{thm: thm1} is straightforward through explicit calculations.

By applying \cite[Theorem~5]{MorrisPeresMR2198701} for $(p^{(n)})^2$, we have the following lemma, which is implied by \cite[Proposition~6.18]{AldousFill2014} for finite regular graphs.
\begin{lem}\label{lem: trace}
Under $(H2)$, we have that for $k\geq 4$,
\[\sharp V^{(n)}\cdot\int |x|^{k}\nu^{(n)}(\mathrm{d}x)\leq 1+4c_1^{-1}c_2^{-4}\cdot\sharp V^{(n)}\cdot \lfloor k/4\rfloor^{-1/2}.\]
In particular, $|\Tr (p^{(n)})^k-1|\leq 4c_1^{-1}c_2^{-4}\cdot\sharp V^{(n)}\cdot \lfloor k/4\rfloor^{-1/2}$.
\end{lem}
\begin{proof}
 Note that $1$ is an eigenvalue of $p^{(n)}$ and that
 \[|\Tr (p^{(n)})^k-1|\leq \sharp V^{(n)}\cdot\int |x|^k\nu^{(n)}(\mathrm{d}x)-1.\]
 Also, note that $\nu^{(n)}$ is supported inside the unit disk centered at $0$ and hence $\int |x|^k\nu^{(n)}(\mathrm{d}x)-1$ is non-increasing in $k$. Hence, it is sufficient to bound $\int |x|^{4k}\nu^{(n)}(\mathrm{d}x)$ for each $k\geq 1$. Since singular values dominate eigenvalues in the $p$-norm (see e.g. \cite[Theorem~3.3.13 b)]{HornJohnsonMR1091716}), we get that $\sharp V^{(n)}\cdot\int |x|^{4k}\nu^{(n)}(\mathrm{d}x)\leq \Tr(A^{t}A)$, where $A^{t}$ is the transpose of $A$ and
 \begin{equation*}
  A=\left[(p^{(n)})^{2k}(x,y)\sqrt{\pi^{(n)}(x)/\pi^{(n)}(y)}\right]_{x,y}.
 \end{equation*}
 (Note that $A$ has the same eigenvalues as $(p^{(n)})^{2k}$.) Following \cite[Eq.~14]{MorrisPeresMR2198701}, for two measures $\nu$ and $\pi$, we write
 \[\chi^2(\nu,\pi)=\sum_{x}(\nu(x))^2/\pi(x)-1.\]
 Then, $\Tr(A^{t}A)$ is exactly $\sum_{x}\pi^{(n)}(x)\chi^2((p^{(n)})^{2k}(x,\cdot),\pi^{(n)})+1$. As we will explain in details in the following, by applying \cite[Theorem~4]{MorrisPeresMR2198701} to $(X^{(n)}_{2j})_{j\geq 0}$, we get that
 \[\max_{x\in V^{(n)}}\chi^2((p^{(n)})^{2k}(x,\cdot),\pi^{(n)})\leq 4c_1^{-1}c_2^{-4}k^{-1/2}\cdot\sharp V^{(n)}\]
 and hence $\sharp V^{(n)}\cdot\int |x|^{4k}\nu^{(n)}(\mathrm{d}x)\leq 1+4c_1^{-1}c_2^{-4}\cdot\sharp V^{(n)}\cdot k^{-1/2}$. To apply \cite[Theorem~4]{MorrisPeresMR2198701}, it suffices to verify that
 \begin{equation}\label{eq: necessary cond for thm4 MorrisPeresMR2198701}
  \int_{4\pi_{*}^{(n)}}^{4/\epsilon}\frac{\mathrm{d}u}{u\psi(u)}\leq k\text{ for }\epsilon=4c_1^{-1}c_2^{-4}k^{-1/2}\cdot\sharp V^{(n)},
 \end{equation}
 where $\pi_{*}^{(n)}=\min_{x\in V^{(n)}}\pi^{(n)}(x)$ and $\psi(u)$ is defined in \cite[Eq. 13]{MorrisPeresMR2198701}. By \cite[Lemma~3]{MorrisPeresMR2198701},
 \[\psi(u)\geq\frac{\gamma^2}{2(1-\gamma)^2}\inf\{\phi_S^2:S\subset V^{(n)}\text{ such that }\sum\limits_{x\in S}\pi^{(n)}(x)\leq \min(u,1/2)\},\]
 where $\gamma=\min(\frac{1}{2},\min_{x}(p^{(n)})^2(x,x))$ and
 \[\phi_S =\sum_{x\in S,y\in S^c}\pi^{(n)}(x)(p^{(n)})^2(x,y)/\sum_{z\in S}\pi^{(n)}(z).\]
 Under our assumption (H2), $\gamma\geq c_2^2/2$, $\pi^{(n)}_{*}\geq c_1/\sharp V^{(n)}$ and $\phi_S\geq c_1c_2^2(u\sharp V^{(n)})^{-1}$ for $S\subset V^{(n)}$ with $\sum_{x}\pi^{(n)}(x)\leq u$. (By a slightly more careful estimate, one could get that $\gamma\geq c_2/2$ and $\phi_S\geq c_1c_2/(2u\sharp V^{(n)})$.) Hence, we have that
 \[\psi(u)\geq c_1^2c_2^8(\min(u,1/2)\cdot\sharp V^{(n)})^{-2}/2\geq c_1^2c_2^8u^{-2}(\sharp V^{(n)})^{-2}/2\]
 and \eqref{eq: necessary cond for thm4 MorrisPeresMR2198701} follows.
\end{proof}

Next, we prove Theorem~\ref{thm: thm1} by Lemma~\ref{lem: trace}:
\begin{proof}[Proof of Theorem~\ref{thm: thm1}]
 \begin{itemize}
  \item[a)] By Lemma~\ref{lem: trace} and (H1), we have that
 \[\int |x|^k\nu(\mathrm{d}x)\leq 4c_1^{-1}c_2^{-4}\cdot \lfloor k/4\rfloor^{-1/2}.\]
 Hence, we have that
  \[\int|\log(1-x)|\nu(\mathrm{d}x)\leq 3+\sum_{k\geq 4}k^{-1}\int|x|^{k}\nu(\mathrm{d}x)\leq 3+8\sqrt{2}c_1^{-1}c_2^{-4}.\]
  \item[b)] Note that by \eqref{eq: based loop measure}, $||\mu^{(n)}||=\sum_{k\geq 1}k^{-1}(1+c^{(n)})^{-k}\Tr (p^{(n)})^k$. By (H1), for all fixed $N\geq 1$, we have that \[\lim_{n\rightarrow\infty}(\sharp V^{(n)})^{-1}\sum_{k=1}^{N}k^{-1}(1+c^{(n)})^{-k}\Tr (p^{(n)})^k=\sum_{k=1}^{N}k^{-1}\int x^{k}\nu(\mathrm{d}x).\]
  By Lemma~\ref{lem: trace}, $\exists C=C(c_1,c_2)<\infty$ such that
  \[\lim_{n\rightarrow\infty}(\sharp V^{(n)})^{-1}\sum_{k>N}k^{-1}(1+c^{(n)})^{-k}\left|\Tr (p^{(n)})^k-1\right|\leq CN^{-1/2}.\]
  Since we take $c^{(n)}=e^{-a\cdot\sharp V^{(n)}+o(\sharp V^{(n)})}$, we have that
  \[\lim_{n\rightarrow\infty}(\sharp V^{(n)})^{-1}\sum_{k>N}k^{-1}(1+c^{(n)})^{-k}=a.\]
  Then, Theorem~\ref{thm: thm1} b) follows.
  \item[c)] We give a brief indication and left the details to the reader. By $b)$, it suffices to show that $\mu^{(n)}(|\ell|\geq q_n)\overset{n\rightarrow\infty}{\sim} a\cdot\sharp V^{(n)}$. By the definition of $\mu^{(n)}$, we write that
  \[\mu^{(n)}(|\ell|\geq q_n)=\sum_{k\geq q_n}k^{-1}(\Tr(p^{(n)})^k-1)(1+c^{(n)})^{-k}+\sum_{k\geq q_n}k^{-1}(1+c^{(n)})^{-k}.\]
  The second summand is the major term, which is asymptotically equivalent to $a\cdot\sharp V^{(n)}$ by our assumptions on $q_n$ and $c^{(n)}$. The first term is $o(a\cdot\sharp V^{(n)})$ as $n\rightarrow\infty$ by Lemma~\ref{lem: trace}.
  \item[d),e)] The proofs are similar and we left them to the reader.
 \end{itemize}
\end{proof}
\begin{rem}\label{rem: connection with st}
 By a direct calculation, the total mass of $\mu^{(n)}$ is $-\log\det(I-(1+c^{(n)})^{-1}p^{(n)})$, see \cite[Eq. (2.5)]{LeJanMR2815763}. Note that $\det((1+c^{(n)})\cdot I-p^{(n)})$ is the partition function of weighted spanning trees rooted at the cemetery point $\partial$ of a Markov chain with transition probabilities $(1+c^{(n)})^{-1}p^{(n)}$, see \cite[Section 8.2]{LeJanMR2815763}. The weight of a tree $T$ is the product of weights $p^{(n)}(x,y)$ on edges $(x,y)$ directed to the root with the convention that $p^{(n)}(x,\partial)=c^{(n)}$. Hence, by using the crude lower bound $c_2^{\# V^{(n)}-1}c^{(n)}$ (given by a single tree) of the partition function of weighted spanning trees, we get that $||\mu^{(n)}||\leq \sharp V^{(n)}(\log(1+c^{(n)})-\log c_2)-\log c^{(n)}$ and hence $\int -\log(1-x)\nu(\mathrm{d}x)\leq -\log c_2$. For a reversible chain (i.e. $\pi^{(n)}(x)p^{(n)}(x,y)=\pi^{(n)}(y)p^{(n)}(y,x)$, $\forall x,y\in V^{(n)}$), $\nu^{(n)}$ and $\nu$ are supported on $[-1,1]$ and Theorem~\ref{thm: thm1} a) is improved:
 \[\int |\log(1-x)|\nu(\mathrm{d}x)=\int -\log(1-x)\nu(\mathrm{d}x)\leq \log (2/c_2).\]
\end{rem}

\section{Example: discrete tori}\label{sec: torus}
We calculate the limiting probability measure $\nu$ and $\int -\log(1-x)\nu(\mathrm{d}x)$ for simple random walks $X^{(d,n)}$ on discrete tori $V^{(d,n)}=\mathbb{Z}^d/n\mathbb{Z}^d$. We denote by $p^{(d,n)}$ the corresponding transition probabilities. Then, (H2) holds with $c_1=1$ and $c_2=(2d)^{-1}$. The eigenvalues of $p^{(1,n)}$ are $\cos(\frac{2\pi}{n}),\cos(\frac{4\pi}{n}),\ldots,\cos(\frac{2\pi n}{n})$, see e.g. \cite[Subsection~12.3.1]{LevinPeresElizabethMR2466937}. Since $p^{(d,n)}$ is a product of $p^{(1,n)}$ in the sense \cite[Eq. 12.19]{LevinPeresElizabethMR2466937}, its eigenvalues are $\frac{1}{d}\left(\cos\left(\frac{2\pi}{n}p_1\right)+\cdots+\cos\left(\frac{2\pi}{n}p_d\right)\right)$ where $p_1,\ldots,p_d$ take values in $\{0,\ldots,n-1\}$. Rewrite these real eigenvalues of $p^{(d,n)}$ in non-decreasing order $\lambda^{(d)}_1\leq \ldots\leq \lambda^{(d)}_{n^d}$. Define $\nu_n^d=\frac{1}{n^d}\sum_{i=1}^{n^d}\delta_{\lambda^{(d)}_i}$ and $\tilde{\nu}_n^d=\frac{1}{n^d}\sum_{i=1}^{n^d}\delta_{d\cdot\lambda^{(d)}_i}$. Then $\tilde{\nu}_n^d=(\nu_n^1)^{*d}$. For all $f\in C([-1,1])$,
\[\lim\limits_{n\rightarrow\infty}\int_{-1}^{1} f(x)\nu_n^1(\mathrm{d}x)=\lim\limits_{n\rightarrow\infty}\frac{1}{n}\sum\limits_{p_1=0}^{n-1}f(\cos(2p_1\pi/n))=\int_{0}^{1}f(\cos(2\pi x))\,\mathrm{d}x=\int_{-1}^{1}\frac{f(x)}{\pi\sqrt{1-x^2}}\mathrm{d}x.\]
Hence, for $d\geq 1$, $\nu$ is the convolution of arcsine distributions on $[-1/d,1/d]$:
\[\nu=\lim\limits_{n\rightarrow\infty}\nu_n^d=m^{*d}\text{ where }m(\mathrm{d}x)=1_{x\in[-1/d,1/d]}d\pi^{-1}(1-d^2 x^2)^{-1/2}\mathrm{d}x,\]
\[\int x^{j}\nu(\mathrm{d}x)=1_{j\text{ is even}}\cdot\sum_{j_1,\ldots,j_d\geq 0:j_1+\cdots+j_d=j/2}(2d)^{-j}\frac{j!}{(j_1!)^2\cdots (j_d!)^2},\]
and that
\[\int-\log(1-x)\,\nu(\mathrm{d}x)=\int_{[0,1]^d}-\log\left(1-\frac{1}{d}\sum_{i=1}^{d}\cos(2\pi x_i)\right)\,\mathrm{d}x_1\cdots\,\mathrm{d}x_d,\]
which equals $\log 2$ when $d=1$ and equals $-4G/\pi+2\log 2$ when $d=2$, where $G$ is the Catalan's constant. We refer to \cite{Kasteleyn19611209} or \cite{MontrollMR0174486} for the evaluation when $d=2$. By Theorem~\ref{thm: thm1} and previous calculations, for $d=1$, we have that $\lim_{n\rightarrow\infty}\mu^{(n)}(\log(|\ell|)/n\in\mathrm{d}x)/n=\log 2\cdot\delta_0(\mathrm{d}x)+1_{0<x<a}\cdot\mathrm{d}x$ and that $\lim_{n\rightarrow\infty}\mu^{(n)}(|\ell|=j)/n=\int x^j\nu(\mathrm{d}x)=1_{j\text{ is even}}\cdot 2^{-j}\binom{j}{j/2}$. Hence, for $d=1$, if we take a Poisson point process $\mathcal{L}^{(n)}$ of loops of the intensity $\mu^{(n)}/n$, then as $n\rightarrow\infty$, $\sum_{\ell\in\mathcal{L}^{(n)}}\delta_{\log|\ell|/n}$ converges to a Poisson point process on $[0,a]$ with intensity measure $\log 2\cdot\delta_{0}+\Leb(\mathrm{d}x)$, where $\Leb$ denotes the Lebesgue measure, and $\sum_{\ell\in\mathcal{L}^{(n)}:|\ell|\leq n^{100}}\delta_{|\ell|}$ converges to a Poisson point process on $2\mathbb{Z}_{+}$ with intensity measure $\sum_{j}2^{-j}\binom{2j}{j}\delta_j$. (Indeed, the sequence $(n^{100})_n$ could be replaced by any sequence of positive integers $(k_n)_n$ such that $\lim_{n\rightarrow\infty}k_n=\infty$ and that $\lim_{n\rightarrow\infty}\log k_n/n=0$.)

\section{Example: balls in a regular tree}\label{sec: balls in regular trees}
\begin{figure}
\includegraphics{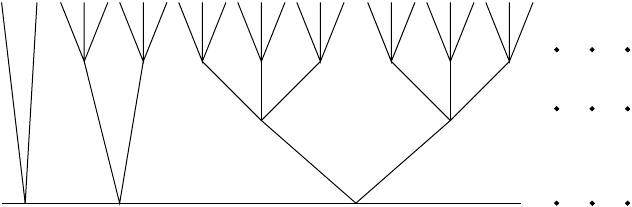}
\caption{Illustration of $\mathcal{T}^{4}$}
\label{fig: canopy}
\end{figure}

We consider the SRW on a truncated regular tree and we will calculate $\nu(\mathrm{d}x)$. More precisely, let $T_d$ be an infinite regular tree with degree $d\geq 3$. Fix a vertex $v_0\in T_d$, let $G^{(d,n)}=B(v_0,n)$ be the balls with radius $n$ centered at $v_0$ and $X^{(d,n)}$ be the SRW on it with transition probabilities $p^{(d,n)}$. Let $\nu^{(d,n)}$ be the empirical distribution of the eigenvalues of $p^{(d,n)}$. By choosing a root vertex $o$ uniformly within $G^{(d,n)}$, we obtain a sequence of random graphs $(G^{(d,n)},o)$, which converges locally to a canopy tree $(\mathcal{T}^{(d)},O)$, where the root $O$ is of distance $k$ from the boundary\footnote{The boundaries are the leaves, i.e. the vertices with degree $1$.} $\partial \mathcal{T}$ with probability $(d-2)(d-1)^{-(1+k)}$ for $k\geq 0$, see Figure~\ref{fig: canopy} for an illustration of $\mathcal{T}^4$. Indeed, by explicit calculations, we have that $d(o,\partial G^{d,n})$ converge in distribution to $d(O,\partial \mathcal{T})$. Moreover, for all positive integers $k$, conditionally on that $o$ is of distance $k$ from the boundary $\partial G^{d,n}$, the distribution of the ball $B(o,k)$ in $G^{d,n}$ equals the ball $B(O,k)$ inside $\mathcal{T}$ for large enough $n$.

For simplicity of notation, we write $\mathcal{T}$ instead of $\mathcal{T}^{(d)}$. The notion of local convergence is introduced in \cite{BenjaminiSchrammMR1873300} and \cite{AldousLyonsMR2354165}. By \cite{BScvgtospeccvg}, it implies the convergence of the spectral measures $\nu^{(d,n)}$. More precisely, $\lim_{n\rightarrow\infty}\nu^{(d,n)}=\nu^{(d)}$ such that
\[\nu^{(d)}=\mathbb{P}[O\text{ is of distance }j\text{ from }\partial\mathcal{T}]\cdot\nu^{(d)}_j,\]
and
$\nu^{(d)}_j(\mathrm{d}z)$ is the probability measure supported on the unit disk such that for vertices $x$ of distance $j$ from $\partial\mathcal{T}$,
\[\int z^k\nu^{(d)}_{j}(\mathrm{d}z)=\mathbb{P}^x[X^{(d)}_k=x]=\lim_{n\rightarrow\infty}\sum_{\substack{y\in G^{(d,n)}\text{ is of}\\
\text{distance }k\text{ from }\partial\mathcal{T}}}\mathbb{P}^{y}[X^{(d,n)}_k=y]/\sharp G^{(d,n)},\]
where $(X^{(d)}_k)_k$ is the SRW on the canopy tree. Hence, $\int z^{k} \nu^{(d)}(\mathrm{d}z)=\mathbb{E}[\mathbb{P}^{O}[X^{(d)}_k=O]]$.

By \cite[Theorem~1.4]{AizenmanWarzelMR2329431}, the transition kernel of the SRW on canopy tree $\mathcal{T}$ acting on $\ell^2(\mathcal{T})$ has only point spectrum with compactly supported eigenfunctions. In that proof, M. Aizenman and S. Warzel used the idea in \cite{AllardFroeseMR1804866} on the decomposition of $\ell^2(\mathcal{T})$ into invariant subspaces. To be more precise, let $(p^{(d)}(x,y))_{x,y\in\mathcal{T}}$ be the transition probabilities of SRW on canopy tree $\mathcal{T}$ and $m$ be a reference measure defined by $m(\{x\})=\deg(x)$ for a vertex $x$ in the canopy tree. Note that $p^{(d)}$ is reversible with respect to $m$ and $\ell^2(\mathcal{T})$ equals $L^2(\mathcal{T},m)$ in the sense of set. As in \cite{AizenmanWarzelMR2329431}, for a vertex $x\in \mathcal{T}$, we define a finite subtree $\mathcal{T}(x)$ at $x$:
\[\mathcal{T}(x)=\{y\in\mathcal{T}: y\text{ is forward (in the direction of }\partial\mathcal{T}\text{) or equal to }x\}\]
and we have an orthogonal decomposition of $L^2(\mathcal{T},m)$ as follows:
\[L^2(\mathcal{T},m)=\bigoplus_{x\in\mathcal{T}}\mathcal{Q}_x\text{ with }\mathcal{Q}_x=\bigoplus_{y\in\mathcal{T}(x),y\sim x}\mathcal{S}_{y}\ominus\mathcal{S}_{x},\]
where $\ominus$ means the orthogonal complement\footnote{To be more precise, $\bigoplus_{y\in\mathcal{T}(x),y\sim x}\mathcal{S}_{y}\ominus\mathcal{S}_{x}=\{f\in \bigoplus_{y\in\mathcal{T}(x),y\sim x}\mathcal{S}_{y}:f\perp \mathcal{S}_x\}$.} and $\mathcal{S}_x$ denotes the subspace of symmetric functions supported on the forward subtree $\mathcal{T}(x)$:
\[\mathcal{S}_{x}\overset{\text{def}}=\{\psi\in L^2(\mathcal{T},m):\psi\text{ is supported on }\mathcal{T}(x)\text{ and is constant on each generation of }\mathcal{T}(x)\}.\]
Then, for each $x\in\mathcal{T}$, $\mathcal{Q}_x$ is an invariant space for the map $p^{(d)}=(p^{(d)}(x,y))_{x,y\in\mathcal{T}}:\psi\rightarrow p^{(d)}\psi$ such that $p^{(d)}\psi(x)=\sum_{y}p^{(d)}(x,y)\psi(y)$. We will describe the eigenvalues and the corresponding eigenfunctions in each $\mathcal{Q}_x$ where $x$ is of distance $N\geq 1$ from $\partial\mathcal{T}$. Consider the following transition probabilities on $\{0,\ldots,N-1\}$:
\[Q_N=\begin{bmatrix}
0 & 1 & 0 & \ldots & 0\\
\frac{d-1}{d} & 0 & \frac{1}{d} & \ddots & \vdots\\
0 & \ddots & \ddots & \ddots & 0\\
\vdots & \ddots & \frac{d-1}{d} & \ddots & \frac{1}{d}\\
0 & \ldots & 0 & \frac{d-1}{d} & 0
\end{bmatrix}.\]
Note that
\begin{equation}\label{eq: det Q by Chebyshev polynomials}
 \det(\lambda\cdot I_N-Q_N)=\left(\frac{\sqrt{d-1}}{d}\right)^N\left(dU_N\left(\frac{d\lambda}{2\sqrt{d-1}}\right)-\frac{d\lambda}{\sqrt{d-1}}U_{N-1}\left(\frac{d\lambda}{2\sqrt{d-1}}\right)\right),
\end{equation}
where $U_N$ are Chebyshev polynomials of the second kind, defined by the identity \[U_N(\cos(\theta))=\sin((N+1)\theta)/\sin(\theta),\quad\forall\theta\in\mathbb{R}.\]
Hence, $Q_N$ has $N$ different eigenvalues $(\lambda_{N,i})_{i=1,\ldots,N}$ which are exactly $2\sqrt{d-1}/d$ times the zeros of
\begin{equation}\label{eq: eig cheby}
 dU_N(x)=2(d-1)xU_{N-1}(x).
\end{equation}
 Denote by $(\varphi_{N,i})_{i=1,\ldots,N}$ the corresponding eigenfunctions. Next, take a $(d-1)\times (d-1)$ real orthogonal matrix $(\xi_{ij})_{i,j=1,\ldots,d-1}$ with $\xi_{1j}=1/\sqrt{d-1}$ for $j=1,\ldots,d-1$. List the neighbor vertices of $x$ in $\mathcal{T}(x)$: $y_1,\ldots,y_{d-1}$. For $N\geq 1$, $i=1,\ldots,N$ and $j=2,\ldots,d-1$, we define a function $\psi_{N,i,j}$ supported on $\cup_{i=1}^{d-1}\mathcal{T}(y_i)$ by taking the value $\xi_{js}\cdot\varphi_{N,i}(k)$ on each $\{z\in\mathcal{T}(y_s): z\text{ is of distance }k\text{ from }\partial\mathcal{T}\}$. Then, $(\psi_{N,i,j})_{j=2,\ldots,d-1}$ are eigenfunctions associated with the eigenvalue $\lambda_{N,i}$ and $(\psi_{N,i,j})_{i=1,\ldots,N;j=2,\ldots,d-1}$ is an orthogonal basis of $\mathcal{Q}_x$. From the spectral representation of the transition probabilities $p^{(d)}$, we obtain the following
\begin{rem}\label{rem: lim canopy}
Let $(\lambda_{N,i})_{i=1,\ldots,N}$ be the zeros of \eqref{eq: eig cheby}. Then,
\[\nu=\sum_{N=1}^{\infty}\frac{(d-2)^2}{(d-1)^{1+N}}\sum\limits_{i=1}^{N}\delta_{\lambda_{N,i}}\]
and $\int-\log(1-x)\nu(\mathrm{d}x)=\sum_{N=1}^{\infty}-\frac{(d-2)^2}{(d-1)^{1+N}}\log\det(I_N-Q_N)$. \footnote{Similarly, by \eqref{eq: det Q by Chebyshev polynomials}, we could express $\int-\log(\lambda-x)\nu(\mathrm{d}x)$ by an infinite sum of functions involving Chebyshev polynomials and logarithms. However, we have no closed form expression for the moments of $\nu$.} Note that $\det(I_N-Q_N)=d^{1-N}$. A quick way to view this is through the connection with random rooted spanning tree with all edges directed towards the root. Indeed, $\det(I-Q_N)$ is the total mass of directed spanning trees rooted at the cemetery point, where the weight of a tree is given by the product of weights $(Q_N)_{xy}$ on directed edges $xy$ in that tree. In this particular case, there is only one rooted tree with weight $d^{1-N}$. Hence, we get that $\int-\log(1-x)\nu(\mathrm{d}x)=(d-1)^{-1}\log d$. If we are only interested in the total mass $||\mu^{(n)}||$, then a simpler way is to use the relation with the spanning trees. Indeed, by \eqref{eq: based loop measure}, we have that
\[||\mu^{(n)}||=-\log\det(I-(1+c^{(n)})^{-1}p^{(d,n)})=\log(1+c^{(n)})\cdot\sharp V^{(n)}+\log\det((1+c^{(n)})\cdot I-p^{(d,n)}),\]
where $V^{(n)}$ is the vertex set of the $n$-th trees in the sequence and $\log\det((1+c^{(n)})\cdot I-p^{(d,n)})$ equals the total mass of directed spanning trees rooted at the cemetery point. (Here, we view $c^{(n)}$ as the killing rate, i.e. the jumping rate to the cemetery point.) For our choice of $c^{(n)}$, the mass is concentrated on trees such that the root has only one neighbor. The total mass of such trees is simply $c^{(n)}\sum_{x\in V^{(n)}}\deg(x)/\prod_{x\in V^{(n)}}\deg(x)$. Hence,
\begin{equation}\label{eq: total mass general tree}
 \lim_{n\rightarrow\infty}||\mu^{(n)}||/\sharp V^{(n)}=a+\sum_{x\in V^{(n)}}\log\deg(x)/\sharp  V^{(n)},
\end{equation}
which equals $a+(\log d)/(d-1)$ for a sequence of balls in a $d$-regular tree. For a sequence of trees in general, when the average degrees are bounded, $\sum_{x\in V^{(n)}}\log\deg(x)/\sharp V^{(n)}$ are also bounded by concavity and \eqref{eq: total mass general tree} holds when $\lim_{n\rightarrow\infty}\sum_{x\in V^{(n)}}\log\deg(x)/\sharp V^{(n)}$ exists. In this case, the uniformly chosen balls $B(o^{(n)},r)\subset G^{(n)}$ with radius $r$ are tight for each $r$ since $\sharp B(o^{(n)},r)$ are tight. (One could show this by induction on $r$ and by using the fact that a uniform chosen neighbor of a uniformly chosen vertex is a uniformly chosen vertex.) Hence, $(G^{(n)})_n$ are tight and locally convergent subsequence exists. One could get similar results along that convergent subsequence.
\end{rem}

\section{Complete graphs}
Let $G^{(n)}$ be the complete graph with vertex set $V^{(n)}$ and $X^{(n)}$ be a SRW on it.
\subsection{Proof of Theorem~\ref{thm: thm2}}
By explicit calculation, we have that
\begin{equation}\label{eq: complete trace}
 \Tr (p^{(n)})^{k}=1+(-1)^k(n-1)^{-(k-1)}.
\end{equation}
By \eqref{eq: based loop measure} and \eqref{eq: complete trace}, we have that
\begin{equation}\label{eq: complete total mass}
 ||\mu^{(n)}||=-\log(c^{(n)}/(1+c^{(n)}))-(n-1)\log(1+(1+c^{(n)})^{-1}(n-1)^{-1})\overset{n\rightarrow\infty}\sim a\log n.
\end{equation}
Take $t\in(0,a)$. We have that
\[\bar{\mu}^{(n)}(\log|\ell|/\log n\leq t)=||\mu^{(n)}||^{-1}\sum\limits_{k\leq n^t}k^{-1}(1+c^{(n)})^{-k}\Tr (p^{(n)})^k.\]
Note that both $(1+c^{(n)})^{-k}$ and $\Tr (p^{(n)})^k$ converge to $1$ uniformly for $k\leq n^t$ as $n\rightarrow\infty$. Hence, $\bar{\mu}^{(n)}(\log|\ell|/\log n\leq t)$ converges to $t/a$ as $n\rightarrow\infty$.

\subsection{Proof of Corollary~\ref{cor: cor of thm2}}\label{subsec: proof of cor1 of thm2}
As before, we denote by $(X^{(n)}_k)_{k\geq 0}$ a simple random walk on the complete graph. We denote by $\mathbb{P}^{(n)}_k$ the law of $(X^{(n)}_j)_{j=0,\ldots,k-1}$ with $X^{(n)}_0$ uniformly distributed on the complete graph. We first compare $\mathbb{P}^{(n)}_k(\mathcal{C})$ with $\bar{\mu}^{(n)}(\mathcal{C}\big| |\ell|=k)$:
\begin{lem}\label{lem: equivalence}
Let $\mathbb{P}^{(n)}_k(\mathcal{C})$ be the probability that $(X_1^{(n)},\ldots,X_k^{(n)})$ covers the complete graph. Then, we have that
\[\frac{n-2}{n}\bar{\mu}^{(n)}(\mathcal{C}\big| |\ell|=k) \leq \mathbb{P}^{(n)}_k(\mathcal{C})\leq\frac{n-1}{n-2}\bar{\mu}^{(n)}(\mathcal{C}\big| |\ell|=k).\]
\end{lem}
\begin{proof}
Define $S^{(k)}_1=\{(x_1,\ldots,x_k)\in\{1,\ldots,n\}^k:x_2\neq x_1,x_3\neq x_2,\ldots,x_k\neq x_{k-1}\}$ and $S^{(k)}_2=\{(x_1,\ldots,x_k)\in\{1,\ldots,n\}^k:x_1\neq x_k,x_2\neq x_1,x_3\neq x_2,\ldots,x_k\neq x_{k-1}\}$. Then, $\sharp S_1^{(k)}=n(n-1)^{k-1}$ and $\sharp S_2^{(k)}=(n-1)^k+(-1)^k(n-1)$. Note that $\bar{\dot\mu}^{(n)}(\cdot\big| |\ell|=k)$ is the uniform distribution on $S_2^{(k)}$ and $\mathbb{P}^{(n)}_k$ is the uniform distribution on $S_1^{(k)}$. Hence,
\[\mathbb{P}^{(n)}_k(\mathcal{C},X^{(n)}_{k}\neq X^{(n)}_{1})=\sharp S^{(k)}_2/\sharp S^{(k)}_1\cdot\bar{\mu}^{(n)}(\mathcal{C}\big| |\ell|=k).\]
Note that $\sharp S_2^{(k)}\leq\sharp S_1^{(k)}\leq \frac{n}{n-2}\sharp S_2^{(k)}$. We see that
\begin{equation}\label{eq: cot1}
\frac{n-2}{n}\bar{\mu}^{(n)}(\mathcal{C}\big| |\ell|=k)\leq\mathbb{P}^{(n)}_k(\mathcal{C},X^{(n)}_{k}\neq X^{(n)}_{1})\leq\bar{\mu}^{(n)}(\mathcal{C}\big| |\ell|=k).
\end{equation}
On the other hand, we have that
\begin{multline}\label{eq: cot2}
\mathbb{P}^{(n)}_k(\mathcal{C},X^{(n)}_{k}=X^{(n)}_{1})=\mathbb{P}^{(n)}_{k-1}\left(\mathbb{P}^{(n)}_k(\mathcal{C},X^{(n)}_{k}=X^{(n)}_{1}|\sigma(X^{(n)}_{1},\ldots,X^{(n)}_{k-1}))\right)\\
=\mathbb{E}^{(n)}_{k-1}\left((n-1)^{-1}\mathbb{P}^{(n)}_{k-1}(\mathcal{C})1_{X^{(n)}_{k-1}\neq X^{(n)}_1}\right)\leq (n-1)^{-1}\mathbb{P}^{(n)}_{k}(\mathcal{C}).
\end{multline}
By combining \eqref{eq: cot1} and \eqref{eq: cot2}, we get that
\[\mathbb{P}^{(n)}_k(\mathcal{C})\leq \bar{\mu}^{(n)}(\mathcal{C}\big| |\ell|=k)+(n-1)^{-1}\mathbb{P}^{(n)}_k(\mathcal{C})\leq \frac{n-1}{n-2}\bar{\mu}^{(n)}(\mathcal{C}\big| |\ell|=k).\qedhere\]
\end{proof}
Consider the cover time $C^{(n)}=\min\{k\geq 0:(X^{(n)}_j)_{j=0,\ldots,k}\text{ contains }V^{(n)}\}$. In the complete graph case, it happens to be a coupon collector problem and it is known that $\frac{C^{(n)}-n\log n}{n}$ converges in law to a Gumbel distribution, see Section 6.2 in \cite{AldousFill2014}. It implies that for all $\delta>0$,
\begin{equation}\label{eq: cot3}
 \lim_{n\rightarrow\infty}\mathbb{P}[C^{(n)}\leq n]=0\text{ and }\lim_{n\rightarrow\infty}\mathbb{P}[C^{(n)}\leq n^{1+\delta}]=1.
\end{equation}
Hence, by Theorem~\ref{thm: thm2}, Lemma~\ref{lem: equivalence} and \eqref{eq: cot3}, we see that $\lim\limits_{n\rightarrow\infty}\frac{\mu_n(\mathcal{C})}{\log n}=(a-1)_{+}$ and Corollary~\ref{cor: cor of thm2} follows.

\paragraph{Acknowledgement}
We thank S. Lemaire, C. Raithel and A. Sapozhnikov for careful readings. We thank Y. Le Jan for useful comments and advices. We thank anonymous readers for careful readings, useful suggestions and valuable comments. We thank anonymous referees for careful readings and useful comments.

\bibliographystyle{amsalpha}
\newcommand{\etalchar}[1]{$^{#1}$}
\providecommand{\bysame}{\leavevmode\hbox to3em{\hrulefill}\thinspace}
\providecommand{\MR}{\relax\ifhmode\unskip\space\fi MR }
\providecommand{\MRhref}[2]{%
  \href{http://www.ams.org/mathscinet-getitem?mr=#1}{#2}
}
\providecommand{\href}[2]{#2}

\end{document}